# ROBOT RELIABILITY USING PETRI NETS AND FUZZY LAMBDA-TAU $(\lambda - \tau)$ METHODOLOGY

AJAY KUMAR[1], S. P. SHARMA[2] AND DINESH KUMAR[3]

*Indian Institute of Technology Roorkee, Roorkee 247 667, India*

Email: [1]ajay1dma@iitr.ernet.in, [2]sspprfma@iitr.ernet.in, [3]dinesfme@iitr.ernet.in.

**ABSTRACT**

Robot reliability has become an increasingly important issue in the last few years due to increased application of robots in many industries (like automobile industry) under hazardous and unstructured environment. As the component failure behavior is dependent on configuration and environment, the available information about the constituent component of robots is most of the time imprecise, incomplete, vague and conflicting and so it is very difficult to analyze their behavior and to predict their failure pattern. The reliability analysis of any system provides an understanding about the likelihood of failures occurring in the system/component and the increased insight about its inherent weakness. The objective of this paper is to quantify the uncertainties that makes the decision more realistic, generic and extendable to application domain. In this paper various reliability parameters (such as mean time between failures, expected number of failures, reliability, availability etc.) are computed using Fuzzy Lambda-Tau methodology. Triangular fuzzy numbers are used to represent failure rates and repair times as they allow expert opinion, linguistic variables, operating conditions, uncertainty and imprecision in reliability information, to be incorporated into system model. Petri Nets are used because unlike the fault tree methodology, the use of Petri Nets allows efficient simultaneous generation of minimal cut and path sets.

## 1. INTRODUCTION

The increasing desire to produce more reliable robots has created interest in several tools used in fault-tolerant design. Such tools seek to evaluate the effectiveness of new designs. The extra component needed for fault-tolerant robot design obviously add extra costs and extra possibility of failure [7, 26]. Therefore reliability analysis is needed to give a hard number showing that the benefits of fault-tolerant design is tangible and worth the effort [2]. Unfortunately, the component failure rates used in these calculations are inaccurate and very often dependent on the configuration and environment, and thus known only approximately [7]. Further age, adverse operating conditions and the vagaries of the system affects each unit of the system differently [3-6]. As such the reliability of the system is affected by various factors such as design, manufacturing, installation, commissioning, operation and maintenance. Consequently it may be extremely difficult if not impossible to construct accurate and complete mathematical model for the system in order to access the reliability because of inadequate knowledge about the basic failure events. In literature, different standard approaches are used during the design phase. The most important of them are Failure Mode and Effect Analysis (FMEA), Fault Tree Analysis (FTA), Failure Mode, Effects and Criticality Analysis (FMECA), Fuzzy Fault tree (FFT), Markov Modeling (MM) and Petri Nets (PNs).

From literature, it is observed that the the standard approaches of reliability engineering rely on the probability model, which is often inappropriate for this task [1, 7]. Probability based analysis usually requires more information about the system than is known, such as mean failure rates or failure rate distribution which commonly results in dubious assumption about the original data. However, fuzzy logic provides necessary requirements in handling with imprecise and uncertain information in more consistent and logical manner.

Among the inexact reasoning methods, fuzzy methodology (FM) acts as one of the most viable and effective tool. Recently, fuzzy methodology has been widely applied in fault diagnosis [8, 9], structural reliability [10], software reliability [11, 12], human reliability [13], safety and risk engineering [8, 14, 15] and quality control [16, 17].

In this paper, a framework to analyze the complex behavior of non-redundant robot [7, 26, 27] using fuzzy lambda-tau ($\lambda - \tau$) methodology [18] has been presented. Petri net [20] model of the robot is obtained from its equivalent fault tree model [26]. Different cut sets are obtained using matrix method [19]. Failure rates and repair times for the robot are computed using fuzzy lambda-tau ($\lambda - \tau$) methodology and various reliability parameters are quantified in terms of fuzzy, crisp and defuzzified values.

Similar to fault tree, Petri nets also make use of diagraph to describe cause and effect relationship between conditions and events. PN has a static as well as dynamic part. The static part consists of only three objects: *places, transitions and arrows*. The dynamic part is the *marking*. Both Petri nets and fault tree are used for software reliability analysis [21], analysis of coherent fault trees [22] and fault diagnosis. In the field of reliability, Petri nets have been presented for reliability evaluation [21], safety analysis [23], Markov Modeling [24] and stochastic modeling [25].

Table-1: Basic expressions of $\lambda - \tau$ Methodology

| Gate | $\lambda_{AND}$ | $\tau_{AND}$ | $\lambda_{OR}$ | $\tau_{OR}$ |
|---|---|---|---|---|
| Expressions | $\prod_{j=1}^{n} \lambda_j \left[ \sum_{i=1}^{n} \prod_{\substack{j=1 \\ i \neq j}}^{n} \tau_j \right]$ | $\dfrac{\prod_{i=1}^{n} \tau_i}{\sum_{j=1}^{n} \left[ \prod_{\substack{i=1 \\ i \neq j}}^{n} \tau_i \right]}$ | $\sum_{i=1}^{n} \lambda_i$ | $\dfrac{\sum_{i=1}^{n} \lambda_i \tau_i}{\sum_{i=1}^{n} \lambda_i}$ |

Table-2: Some Reliability parameters

| Parameters | Expressions |
|---|---|
| Mean Time to Failure | $MTTF_s = \dfrac{1}{\lambda_s}$ |
| Mean Time to Repair | $MTTR_s = \dfrac{1}{\mu_s} = \tau_s$ |
| Mean Time Between Failures | $MTBF_s = MTTF_s + MTTR_s$ |
| Availability | $A_s = \dfrac{\mu_s}{\lambda_s + \mu_s} + \dfrac{\lambda_s}{\lambda_s + \mu_s} e^{-(\lambda_s + \mu_s)t}$ |
| Reliability | $R_s = e^{-\lambda_s t}$ |
| Expected Number of Failures | $W_s = \dfrac{\lambda_s \mu_s t}{\lambda_s + \mu_s} + \dfrac{\lambda^2{}_s}{(\lambda_s + \mu_s)^2}[1 - e^{-(\lambda_s + \mu_s)t}]$ |

In this paper, the static part of Petri nets is used to model the qualitative behavior of non-redundant robot [26]. Minimal cut sets for calculating the reliability parameters are obtained using matrix method. The basic expressions for fuzzy Lambda-Tau methodology and some reliability parameters are given in Table [1-2] for ready reference.

## 2. A PRACTICAL CASE

The present work is based on calculating the reliability of non-redundant robot with two joints and one sensor per joint i.e. it actually has two motors ($M_1$ and $M_2$) and two sensors ($S_1$ and $S_2$). The Fault tree model and equivalent Petri net model of the robot is depicted in Figure 1.

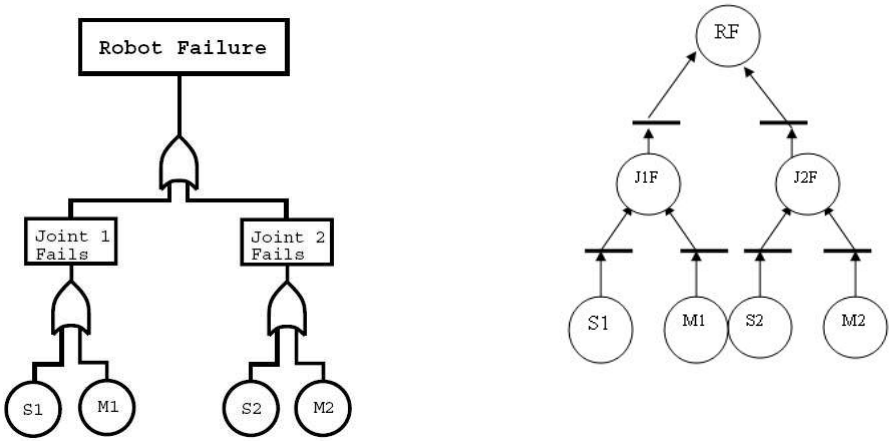

Figure-1:   (a) Fault tree model                                                         (b) Equivalent Petri Net Model

Minimum cut sets are calculated using matrix method and they are $\{S_1\}$, $\{M_1\}$, $\{S_2\}$ and $\{M_2\}$. The following assumptions were taken while modeling the system:

- Component failures and repairs are statistically independent, constant and obey exponential distribution.
- After repairs, the repaired component is considered as good as new.

The procedural steps of the proposed methodology are given below:

**Step 1:** The data related to failure rates $\lambda_i$ and repair times $\tau_i$ of the components $i = 1$ (sensor $S_1$), $i = 2$ (Motor $M_1$), $i = 3$ (Sensor $S_2$) and $i = 4$ (Motor $M_4$) are collected from the historical/present records of the system as presented in Table-3.

Table-3: Failure and repair data

| | |
|---|---|
| $\lambda_1 = \lambda_3 = 0.000182$ failures/h | $\lambda_2 = \lambda_4 = 0.0092$ failures/h |
| $\tau_1 = \tau_3 = 3$ h | $\tau_2 = \tau_4 = 5$ h |

**Step 2:** To handel the vagueness and uncertainty in data, the crisp data of $\lambda$ and $\tau$ is converted into triangular fuzzy number [30] with $\pm 15\%$ ($\pm 25\%$, $\pm 50\%$; as depicted in Figure 2).

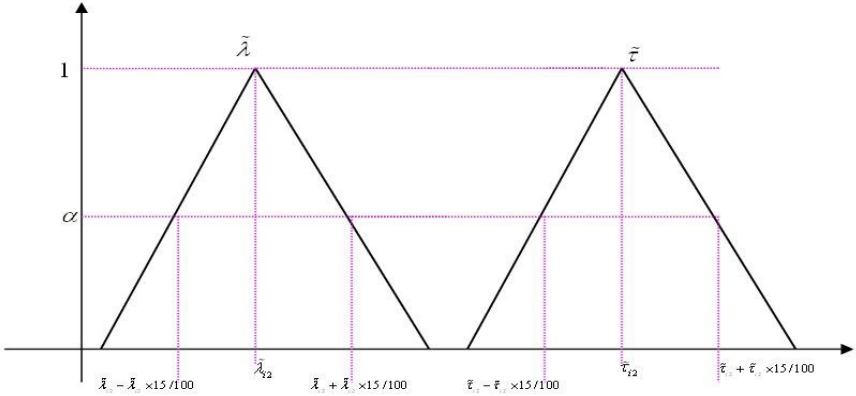

Figure 2: Input Fuzzy Triangular Numbers for failure rate and repair time

**Step 3:** As soon as, the input fuzzy triangular numbers for failure rates and repair times for each of the components are known, the corresponding fuzzy value ($\tilde{\lambda}$ and $\tilde{\tau}$) of crisp failure rate $\lambda$ and repair time $\tau$ can be obtained using extension principle coupled with $\alpha$--cut and interval arithmetic operations on fuzzy triangular numbers [27, 28]. To analyze the system behavior qualitatively, various reliability parameters such as failure rate, repair time, availability, MTBF, reliability and expected number of failures, with left and right spreads are computed at various membership grade and shown graphically in Figure 3.

**Step 4:** It is necessary to convert the fuzzy output to a crisp value as most of the actions or decisions implemented by humans or machines are binary or crisp. The process of converting fuzzy output to a crisp value is said to be *defuzzification*. There exist many defuzzification techniques in the literature [30] such as *max-membership principle, center of area COA, center of sum, center of largest area etc.*, which can be used depending on the application. The COA method is selected for this study as it is equivalent to mean of data and so it is very appropriate for reliability calculation. It the membership function $\mu_{\tilde{A}}(x)$ of the output fuzzy set $\tilde{A}$ is described on the interval $[x_1, x_2]$, then COA defuzzification $\bar{x}$ can be defined as:

$$\bar{x} = \frac{\int_{x_1}^{x_2} x \cdot \mu_{\tilde{A}}(x) dx}{\int_{x_1}^{x_2} \mu_{\tilde{A}}(x) dx}$$

## 3. BEHAVIOR ANALYSIS

The crisp and defuzzified values for various reliability parameters at $\pm 15\%$, $\pm 25\%$ and $\pm 50\%$ spreads are calculated and depicted in Table-4 which clearly indicates that the defuzzified values of various reliability parameters change with change of spreads. It also shows that the defuzzified values of failure rate, repair time and expected number of failures increase with the increase of spreads. On the other hand the defuzzified values of mean timebetween failures, reliability and availability decrease with increase of spreads. Thus from the above analysis it is clear that the maintenance should be based on defuzzified MTBF rather than on crisp MTBF because by defuzzified value of MTBF a safe interval between times of maintenance can be established and inspections can be conducted to monitor the condition or status of various equipments of the system before it reaches the crisp value.

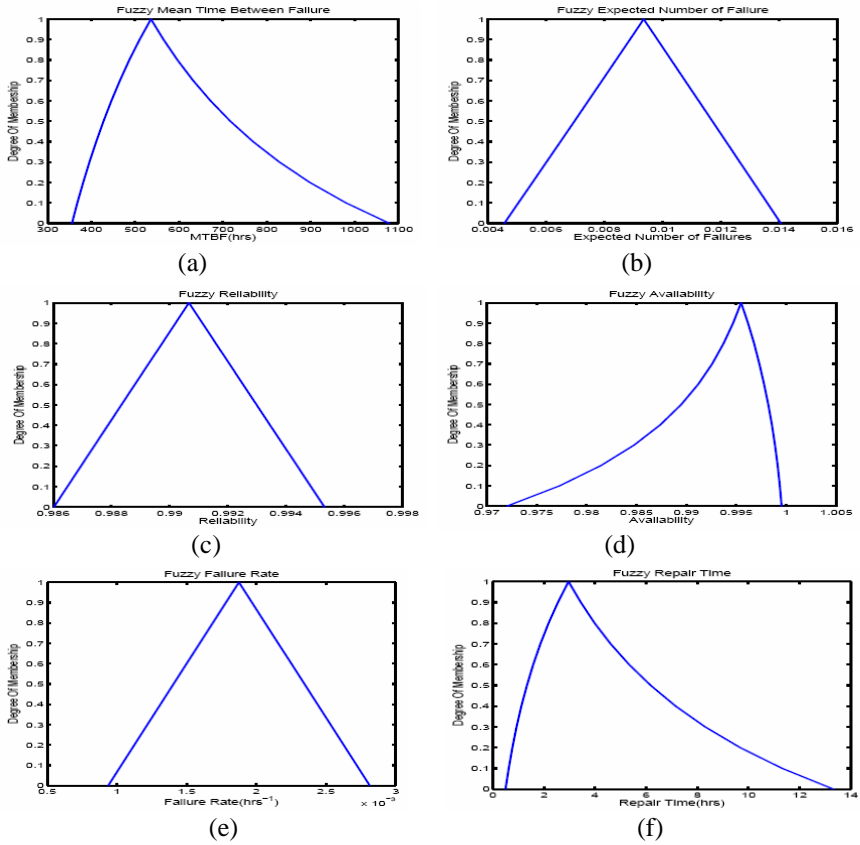

Figure 3: Reliability Parameters for Non-Redundant Robot

Table-4: Crisp and defuzzified values at different spreads

| Reliability Parameters | Crisp values | Defuzzified values ($\pm 15\%$) | Defuzzified values ($\pm 25\%$) | Defuzzified values ($\pm 50\%$) |
|---|---|---|---|---|
| Failure rate ($h^{-1}$) | $1.828412 \times 10^{-3}$ | $1.876400 \times 10^{-3}$ | $1.952427 \times 10^{-3}$ | $2.1315275 \times 10^{-3}$ |
| Repair rate (h) | 2.961202 | 3.097652 | 3.354636 | 4.878753 |
| ENOF | $9.355101 \times 10^{-3}$ | $9.553484 \times 10^{-3}$ | $9.859787 \times 10^{-3}$ | $10.033682 \times 10^{-3}$ |
| MTBF (h) | $5.358966 \times 10^{2}$ | $5.229019 \times 10^{2}$ | $4.939305 \times 10^{2}$ | $4.4216472 \times 10^{2}$ |
| Reliability | 0.990667 | 0.990663 | 0.990662 | 0.990661 |
| Availability | 0.995485 | 0.995219 | 0.994681 | 0.990835 |

## 4. CONCLUSION

In this paper a structured framework has been developed that may help the maintenance engineers to analyze and predict the system behavior. An attempt has also been made to deal with imprecise, uncertain dependent information related to system performance. Various reliability parameters (such as failure rate, repair time, mean time between failures, availability, reliability and expected number of failures) were computed to predict the system behavior in objective terms and it is concluded that in order to improve the availability and reliability aspects, it is necessary to enhance the maintainability requirement of the system.

## 5. ACKNOWLEDGEMENT

One of the authors, **Ajay Kumar** acknowledges the **Council of Scientific and Industrial Research**, **New Delhi**, for providing financial support during the period of preparation of this paper.